\documentclass[11pt]{article}

\usepackage[small,bf]{caption}
\usepackage{amsmath, amsthm, amssymb}
\usepackage{algorithmic}
\usepackage{algorithm}
\usepackage[parfill]{parskip}
\usepackage{soul}

\usepackage{epsfig,fullpage,psfrag,url,verbatim}
\usepackage{graphics}
\usepackage{graphicx}
\usepackage{epstopdf}

\usepackage{amsfonts}
\usepackage{color}

\newtheorem{cor}{Corollary}[section]
\newtheorem{thm}{Theorem}[section]
\newtheorem{lem}{Lemma}[section]
\newtheorem{prop}{Proposition}[section]

\newtheorem{mydef}{Definition}[section]
\newtheorem{rem}{Remark}[section]

\newcommand{\cx}{\check{x}}
\newcommand{\barx}{\bar{x}}

\newcommand{\hh}{h}
\newcommand{\ff}{f}
\newcommand{\EE}{\mathbb{E}}
\newcommand{\PP}{\mathbb{P}}

\newcommand{\hx}{\hat{x}}
\newcommand{\Dh}{D_h}
\newcommand{\intQ}{\mathrm{int}\ Q}
\newcommand{\tg}{\tilde{g}}
\newcommand{\tM}{\tilde{M}}

\newcommand{\cQ}{\mathcal{Q}}

\newcommand{\RR}{\mathbb{R}}

\begin{document}

\title{``Relative-Continuity'' for Non-Lipschitz Non-Smooth Convex Optimization using Stochastic (or Deterministic) Mirror Descent}

\author{Haihao Lu\thanks{MIT Department of Mathematics and MIT Operations Research Center, 77 Massachusetts Avenue, Cambridge, MA   02139
({haihao@mit.edu}).  The author's research is supported by AFOSR Grant No. FA9550-15-1-0276 and the MIT-Belgium Universit\'{e} Catholique de Louvain Fund.}
}
\date{August 13, 2018} 

\maketitle

\begin{abstract}
The usual approach to developing and analyzing first-order methods for non-smooth (stochastic or deterministic) convex optimization assumes that the objective function is uniformly Lipschitz continuous with parameter $M_f$. However, in many settings the non-differentiable convex function $f$ is not uniformly Lipschitz continuous -- for example (i) the classical support vector machine (SVM) problem, (ii) the problem of minimizing the maximum of convex quadratic functions, and even (iii) the univariate setting with $f(x) := \max\{0, x\} + x^2$. Herein we develop a notion of ``relative continuity'' that is determined relative to a user-specified ``reference function'' $h$ (that should be computationally tractable for algorithms), and we show that many non-differentiable convex functions are relatively continuous with respect to a correspondingly fairly-simple reference function $h$. We also similarly develop a notion of ``relative stochastic continuity'' for the stochastic setting. We analyze two standard algorithms -- the (deterministic) mirror descent algorithm and the stochastic mirror descent algorithm -- for solving optimization problems in these new settings, providing the first computational guarantees for instances where the objective function is not uniformly Lipschitz continuous. This paper is a companion paper for non-differentiable convex optimization to the recent paper by Lu, Freund, and Nesterov, which developed analogous results for differentiable convex optimization.
\end{abstract}

\section{Introduction}
The usual  approach to developing and analyzing first-order methods for non-differentiable convex optimization (which we  review shortly) assumes that the objective function is uniformly Lipschitz continuous, in both deterministic and stochastic settings. However, in many settings the objective function $f$ is not uniformly Lipschitz continuous.  For example, consider the Support Vector Machine problem (SVM) for binary classification in machine learning, whose optimization formulation is:
\begin{equation*}
\mbox {SVM:} \ \ \ \min_x f(x):= \tfrac{1}{n} \sum_{i=1}^n \max\{0, 1-y_i x^T w_i\} + \tfrac{\lambda}{2}\|x\|_2^2 \ ,
\end{equation*}

where $w_{i}$ is the input feature vector of sample $i$ and $y_{i} \in \{-1,1\}$ is the label of sample $i$.  Notice that $\ff$ is not differentiable due to the presence of hinge loss terms in the summation, and $\ff$ is also not Lipschitz continuous due to the presence of the $\ell_2$-norm regularization term; thus we cannot directly utilize typical subgradient or gradient schemes and their associated computational guarantees for SVM.  

Another example is the the problem of computing a point $x\in \RR^m$ in the intersection of $n$ ellipsoids, which can be tackled via the optimization problem 
\begin{equation*}
\mbox{IEP:} \ \ \ f^* = \min_x f(x) := \max_{0 \le i \le n} \{ \tfrac{1}{2} x^T A_i x + b_i^T x + c_i \} \ ,
\end{equation*}
where the $i^{\mathrm{th}}$ ellipsoid is $\cQ_i= \{x\in \RR^m : \tfrac{1}{2}x^T A_i x + b_i x + c_i \le 0\}$ and $A_i\in \RR^{m\times m}$ is a symmetric positive semi-definite matrix, $i = 1, \ldots, n$.  Observe that the objective function $\ff$ is both non-differentiable and non-Lipschitz, and so it falls outside of the scope of standard classes of optimization problems for which first-order methods are guaranteed to work. Nevertheless, using the machinery developed in this paper, we will show in Section \ref{wrightwoods} how to solve both of these problems using deterministic or stochastic Mirror Descent.	

In this paper we develop a general theory and algorithmic constructs that overcome the drawbacks in the usual analyses of first-order methods that are grounded on restricted notions of uniform Lipschitz continuity.  Here we develop a notion of ``relative continuity'' with respect to a given convex ``reference function'' $h$, a notion which does not require the specification of any particular norm --  indeed $h$ need not be strongly (or even strictly) convex.  Armed with ``relative continuity'', we demonstrate the capability to solve a more general class of non-differentiable convex optimization problems (without uniform Lipschitz continuity) in both deterministic and stochastic settings.

This paper is a companion for non-differentiable convex optimization to our predecessor paper \cite{lu2018relatively} for differentiable convex optimization. In \cite{lu2018relatively}, with a very similar philosophy, we developed the notion of relative smoothness and relative strong convexity with respect to a given convex reference function.  In that paper we demonstrated the capability to solve a more general class of differentiable convex optimization problems (without uniform Lipschitz continuous gradients), and we also demonstrated linear convergence results for a Primal Gradient Scheme when the objective function $\ff$ is both smooth and strongly convex.

There are some concurrent works on smooth optimization sharing a similar spirit to \cite{lu2018relatively}. Bauschke, Bolte, and Teboulle \cite{bbt} presents a similar definition of relative smoothness as in \cite{lu2018relatively} and analyzes the convergence of Mirror Descent Algorithm, although their algorithm and convergence complexity depend on a symmetry measure of the Bregman distance.  Zhou et al. \cite{zhouaunified} discusses a unified proof of Mirror Descent and the Proximal Point Algorithm under a similar assumption of relative smoothness.  Nguyen \cite{van2017forward} develops similar ideas on analyzing Mirror Descent in a Banach space. A more detailed discussion comparing these related works is also presented in \cite{lu2018relatively}. More recently, Hanzely and Richtarik \cite{hanzely2017randomized} develop stochastic algorithms for the relatively smooth optimization setting.

In Section 2 we review the traditional set-up for Mirror Descent in both the deterministic and stochastic settings. In Section 3 we introduce our notion of ``relative continuity'' in both the deterministic and stochastic settings, together with some relevant properties.  In Section 4 we prove computational guarantees associated with the Mirror Descent and Stochastic Mirror Descent algorithms under relative continuity.  In Section 5 we show constructively how our ideas apply to a large class of non-differentiable and non-Lipschitz convex optimization problems that are not otherwise solvable by traditional first-order methods. Also in Section 5 we analyze computational guarantees associated with Mirror Descent and Stochastic Mirror Descent for the Intersection of Ellipsoids Problem (IEP) and also the Support Vector Machine (SVM) problem. 

\noindent {\bf Notation.} $\|\cdot\|$ denotes a given norm in $\RR^n$ and $\|\cdot\|_*$ denotes the usual dual norm on the dual space. $\| x \|_2 := \sqrt{x^Tx}$ denotes the Euclidean (inner product) norm, where $x^T$ means the transpose of the vector $x$, and $B_2(c,r) := \{ x \in \RR^n : \|x-c\|_2 \le r \}$. $\| A \|_2$ denotes the $\ell_2$ (spectral) norm of a matrix $A$. The inner product $\langle \cdot, \cdot \rangle$ specifically denotes the dot inner product in the underlying vector space. For a conditional random variable $s(x)$ given $x$, $\EE[s(x)|x]$ denotes the conditional expectation of $s(x)$ given $x$.

\section{Traditional Mirror Descent}

The optimization problem of interest is:
\begin{equation}\label{poi1}
\begin{array}{lrlr} P:  &  \ \  \min_x & f(x) \\ \\
& \mbox{ s.t. } & x \in Q \ ,
\end{array}
\end{equation}
where $Q \subseteq \RR^n$ is a closed convex set and $f : Q \to \mathbb{R}$ is a convex function that is not necessarily differentiable.  There are very many deterministic and stochastic first-order methods for tackling \eqref{poi1}, see for example \cite{bubeck2015convex}, \cite{nesterovBook}, \cite{nedic2014stochastic} and the references therein. Virtually all such methods are designed to solve \eqref{poi1} when the objective function $f$ satisfies a uniform Lipschitz continuity condition on $Q$, which in the deterministic setting is (essentially) equivalent to the condition that there exists a constant $M_f < \infty$ for which:
\begin{equation}\label{lipschitz}
\| g(x) \|_* \le M_f  \ \ \ \text{for~all~} x \in Q \ \text{and} \ g(x)\in\partial f(x) \ ,
\end{equation} where $\partial f(x)$ is the subdifferential of $\ff$ at $x$ (i.e., the collection of subgradients of $\ff$ at $x$), $\|\cdot\|$ is a given norm on $\RR^n$, and $\|\cdot\|_*$ denotes the usual dual norm.  

Here we use ``$g(x)$'' to denote an assignment of a subgradient (or an oracle call thereof) at $x$, and so $g(x)$ is not a function nor is it a point-to-set map.  

Another useful functional notion is strong convexity:  $\ff$ is (uniformly) $\mu_f$-strongly convex for some $\mu_f >0$ if
\begin{equation}\label{arm}
f(y) \ge f(x) + \langle g(x), y-x \rangle + \tfrac{\mu_f}{2} \|y-x\|^2 \   \ \ \ \text{for~all~} x,y \in Q \ \text{and} \ g(x)\in\partial f(x) \ .
\end{equation}
Nedi\'{c} and Lee \cite{nedic2014stochastic} obtain improved convergence guarantees for the stochastic mirror descent algorithm under strong convexity, for example.  

\subsection{Deterministic Setting}

Let us now recall the Mirror Decent Algorithm (see \cite{NemirovskyYudin83} and \cite{beckteb03mirror}), which is also referred to as the prox subgradient method when interpreted in the space of primal variables.  Mirror Descent employs a differentiable convex ``prox function'' $h$ to define a Bregman distance:
\begin{equation}\label{eq:bregman}
\Dh(y,x) := h(y) - h(x) - \langle \nabla h(x), y-x \rangle \ \ \ \text{for~all~} x,y \in Q \ . \end{equation}
The Bregman distance is used in the computation of the Mirror Descent update:
\begin{equation*}
x^{i+1} \gets \arg\min_{x \in Q} \left\{f(x^i) + \langle g(x^i), x-x^i \rangle+ \tfrac{1}{t_i}\Dh(x,x^i)\right\} \ ,
\end{equation*} where $\{t_i\}$ is the sequence of step-sizes for the scheme.  A formal statement of the Mirror Descent Algorithm is presented in Algorithm \ref{al:proxgradscheme}.  The traditional set-up requires that $\hh$ is $1$-strongly convex with respect to the given norm $\|\cdot\|$, and in this set-up one can prove that after $k$ iterations it holds for any $x \in Q$ that:
\begin{equation}\label{sd_bound3}
 \min_{0\le i\le k} f(x^i)  - f(x) \ \ \leq \ \  \frac{\tfrac{1}{2}M_f^2\sum_{i=0}^k t_i^2 + \Dh(x,x^0)}{\sum_{i=0}^k t_i} \ ,
\end{equation} which leads to an $O(1/\sqrt{k})$ sublinear rate of convergence using an appropriately chosen step-size sequence $\{t_i\}$, see \cite{beckteb03mirror}.

\begin{algorithm}
\caption{Deterministic Mirror Descent Algorithm with Bregman distance $\Dh(\cdot,\cdot)$}\label{al:proxgradscheme}
$ \ $
\begin{algorithmic}
\STATE {\bf Initialize.}  Initialize with $x^0 \in Q$.  Let  $\hh$ be a given convex function. \\
At iteration $i$ :\\
\STATE  {\bf Perform Updates.}  Compute $g(x^i) \in \partial f(x^i)$ , determine step-size $t_i$, and compute update:
\begin{equation}\label{proxgrad}
x^{i+1} \gets \arg\min_{x \in Q} \{f(x^i) + \langle g(x^i), x-x^i \rangle + \tfrac{1}{t_i} \Dh(x,x^i) \} \ .
\end{equation}
\end{algorithmic}
\end{algorithm}\medskip

Notice in \eqref{proxgrad} by construction that the update requires the capability to solve instances of a ``linearized subproblem'' (which we denote by ``LS'') of the general form:
\begin{equation}\label{subproblem}
\mathrm{LS:} \ \ \ \ \ x_{\text{new}} \gets \arg\min_{x \in Q}\{\langle c, x \rangle + h(x)\} \ ,
\end{equation} for suitable iteration-specific values of $c$.  Indeed, \eqref{proxgrad} is an instance of \eqref{subproblem} with $c = t_i g(x^i) - \nabla h(x^i)$ at iteration $i$.  It is especially important to note that the Mirror Descent update \eqref{proxgrad} is somewhat meaningless absent the capability to efficiently solve \eqref{subproblem}, a fact which we return to later.  In the usual design and implementation of Mirror Descent for solving \eqref{poi1}, one attempts to specify the norm $\|\cdot\|$ and the $1$-strongly convex prox function $h$ in consideration of properties of the feasible domain $Q$ while also ensuring that the LS subproblem \eqref{subproblem} is efficiently solvable.

Notice that the Mirror Descent Algorithm (Algorithm \ref{al:proxgradscheme}) itself does not require the traditional set-up that $\hh$ be $1$-strongly convex for some particular norm; rather this requirement is part of the traditional \emph{analysis}.  As we will see, we can instead analyze Mirror Descent by considering the intrinsic ways that $\ff$ and $\Dh(\cdot,\cdot)$ are related functionally, in a manner that is constructive in terms of actual algorithm design and implementation.  Furthermore, this is in the same spirit as was done in the predecessor paper \cite{lu2018relatively}.

\subsection{Stochastic Setting}
For some convex functions, computing an exact subgradient at $x \in Q$ may be expensive or even intractable, but sampling a random stochastic estimate of a subgradient at $x$, which we denote by $\tg(x)$, may be easy.  We say that $\tg(x)$ is an unbiased stochastic subgradient if $\EE \left[\tg(x) | x\right] \in \partial f(x)$.  The usefulness of a stochastic subgradient methodology is easily seen in the context of machine and statistical learning problems.  A prototypical learning problem is to compute an approximate solution of the following empirical loss minimization problem:
\begin{equation}\label{poi1stoc}
\begin{array}{rlr}  \ \  \min_x & f(x) := \tfrac{1}{n} \sum_{j=1}^n f_j(x) \\ \\
\mbox{ s.t. } & x \in Q \ ,
\end{array}
\end{equation}
where $f_j$ is a non-differentiable convex loss function associated with sample $j$, for $j=1, \ldots, n$ data samples.  When $n \gg 0$, the standard subgradient method needs to evaluate $n$ subgradients in order to compute a subgradient of $\ff$, which can be prohibitively expensive.  A typical alternative is to compute a stochastic subgradient.  Letting $x^i$ denote the $i^{\mathrm{th}}$ iterate, at iteration $i$ a single sample index $\tilde j$ is drawn uniformly and independently on $\{1,\ldots,n\}$, and then a subgradient $\tilde g \in \partial f_{\tilde j}(x^i)$ is computed that is used to define $\tilde g(x^i):= \tilde g$.  This stochastic subgradient is then used in place of a subgradient at iteration $i$.
Notice that by construction $\tg(x^i)$ is a conditional random variable given $x^i$, and $\tg(x^i)$ is an unbiased stochastic subgradient, namely $\EE [\tg(x^i)|x^i] \in \partial f(x^i)$.  

A stochastic version of Mirror Descent is presented in Algorithm \ref{al:proxstocgradscheme}.  The structure of Stochastic Mirror Descent is identical to that of Mirror Descent, the only difference being that the stochastic estimate of a subgradient $\tilde g(x^i)$ replaces the exact subgradient $g(x^i)$ in Algorithm \ref{al:proxstocgradscheme}.

\begin{algorithm}
\caption{Stochastic Mirror Descent Algorithm with Bregman distance $\Dh(\cdot, \cdot)$}\label{al:proxstocgradscheme}
$ \ $
\begin{algorithmic}
\STATE {\bf Initialize.}  Initialize with $x^0 \in Q$.  Let $\hh$ be a given convex differentiable function. \\
At iteration $i$ :\\
\STATE  {\bf Perform Updates.}  Compute a stochastic subgradient $\tg(x^i)$, determine step-size $t_i$, and compute update:\\\medskip
\ \ \ \ \ \ \ \ \ \ \ \ \ \ \ \ \ \ \ \ \ \ \ \ \ \ \ \ $x^{i+1} \gets \arg\min_{x \in Q} \{f(x^i) + \langle \tg(x^i), x-x^i \rangle + \tfrac{1}{t_i} \Dh(x,x^i) \}$ .\\\medskip
\end{algorithmic}
\end{algorithm}\medskip

A standard condition that is required in the traditional convergence analysis for Stochastic Mirror Descent (as well as other  stochastic first-order methods)  is that there exists $G_f >0$ for which:
\begin{equation}\label{def:tra_stoc_cont}
\EE [\|\tg(x)\|_*^2|x] \le G_f^2, \  \text{for any} \ x\in Q \ .
\end{equation}
For notational convenience, we say that $\ff$ is $G_f$-stochastically continuous if \eqref{def:tra_stoc_cont} holds. In \cite{nedic2014stochastic}, Nedi\'{c} and Lee developed convergence results for Stochastic Mirror Descent (Algorithm 2). Under the conditions that (i) $\ff$ is $G_f$-stochastically continuous, (ii) $\hh$ is a differentiable and $\mu_h$-strongly convex function on $Q$, and (iii) $Q$ is a closed bounded set, Nedi\'{c} and Lee (\cite{nedic2014stochastic} equation (27)) show the following convergence result using step-sizes $t_i = \sqrt{\frac{\mu_h D_{\max}}{G_f(i+1)}}$ :
\begin{equation}\label{sto_st_nonstrong_bound}
 \EE\left[  f(\barx^k)\right]  - f^* \ \ \leq \ \  \frac{3G_f \sqrt{ D_{\max}}}{2\sqrt{\mu_h(k+1)}} \ ,
\end{equation}
where $\barx^k:=\frac{1}{\sum_{i=0}^k t_i} \sum_{i=0}^kt_i x^i$ and $D_{\max}:= \max_{x,y\in Q} D_h(x,y)$.

Furthermore, if also (a) $\ff$ is $\mu_f$-strongly convex, and (b) $\hh$ is $L_h$-smooth, Nedi\'{c} and Lee (\cite{nedic2014stochastic} Theorem 1) show that with step-sizes $t_i=\tfrac{2L_h}{\mu_f (i+1)}$ it holds that:
\begin{equation}\label{sto_st_strong_bound}
 \EE\left[ f(\cx^k) \right]  - f^* \ \ \leq \ \  \frac{2G_f^2L_h}{\mu_f(k+1)\mu_h} \ ,
\end{equation}
where $\cx^k:=\frac{2}{(i+1)(i+2)}\sum_{i=0}^k(i+1)x^i$.

\section{Relative Continuity}

In this section we introduce our definition of relative continuity of a function $\ff$ -- actually two different definitions -- one for the deterministic and another for the stochastic setting.  The starting point is a ``reference function'' $\hh$ which is a given differentiable convex function on $Q$ that is used to construct the usual Bregman distance $\Dh(\cdot,\cdot)$ \eqref{eq:bregman}, and that is used as part of the Mirror Descent update \eqref{proxgrad}.  However, we point out for emphasis that unlike the traditional set-up there are no assumptionss on $\hh$ (such as strong or strict convexity).

\subsection{Deterministic Setting}

Consider the objective function $\ff$ of \eqref{poi1}.  We define ``relative continuity'' of $\ff$ relative to the reference function $\hh$ using the Bregman distance $\Dh(\cdot,\cdot)$ of $\hh$ as follows.\medskip

\begin{mydef}\label{def:continuity}
$\ff$ is $M$-relative continuous with respect to the reference function $\hh$ on $Q$ if for any $x,y\in Q$, $x \ne y$, and $g(x)\in \partial f(x)$, it holds that
\begin{equation}\label{eq:cont}
\|g(x)\|_* \le \frac{M\sqrt{2D_h(y,x)}}{\|y-x\|} \ .
\end{equation}

\end{mydef}

(In the particular case when $h(x) = \tfrac{1}{2} \|x\|_2^2$, the Bregman distance is $D_h(y,x) = \tfrac{1}{2} \|y-x\|_2^2$, and the relative continuity condition \eqref{eq:cont} becomes $ \|g(x)\|_2 \le M $, which corresponds to the standard definition of Lipschitz continuity \eqref{lipschitz} for the $\ell_2$-norm.)

We can rewrite \eqref{eq:cont} as
\begin{equation}\label{eq:contalt}
\|g(x)\|^2_* \le M^2 \frac{D_h(y,x)}{\frac{1}{2}\|y-x\|^2} \ ,
\end{equation}
which states that the square of the norm of any subgradient is bounded by the ratio of the Bregman distance $D_h(y,x)$ to $\tfrac{1}{2}\|y-x\|^2 $. \medskip

The following proposition presents the ``key property'' of an $M$-relative continuous function that is used in the proofs of results to follow.\medskip

\begin{prop}\label{prop:equiv}
{\bf (Key property of $M$-relative continuity)} If $\ff$ is $M$-relative continuous with respect to the reference function $\hh$, then for any $t>0$ it holds for all $x,y \in Q$ and $g(x) \in \partial f(x)$ that:
\begin{equation}\label{eq:prop31}
\tfrac{1}{t}D_h(y,x) + \langle g(x), y-x \rangle + \tfrac{1}{2} tM^2 \ge 0 \ .
\end{equation}

\end{prop}

{\bf Proof:}
If $\ff$ is $M$-relative continuous with respect to $\hh$, then for any $t>0$ it follows that
\begin{equation*}
- \langle g(x), y-x \rangle\le \|g(x)\|_* \|y-x\| \le M \sqrt{2D_h(y,x)} \le \tfrac{1}{2} tM^2 + \frac{D_h(y,x)}{t} \ ,
\end{equation*}
where the last inequality is an application of the arithmetic-geometric mean inequality. The proof follows by rearranging terms.  \qed

The ``key property'' \eqref{eq:prop31} is what is used in the proofs of results herein, so we could define $M$-relative continuity using \eqref{eq:prop31} instead of \eqref{eq:cont}.  Furthermore, \eqref{eq:prop31} is independent of any norm structure, and so is attractive for its generality.  However, we use the definition \eqref{eq:cont} because it leads to easy verification of $M$-relative continuity in practical instances as is shown in Section \ref{wellfleet}.



The following proposition presents some scaling and additivity properties of relative continuity.\medskip

\begin{prop}\label{prop:additivity}{\bf Additivity of Relative Continuity}\
\begin{description}
\item[1.]If $\ff$ is $M$-relative continuous with respect to $\hh$, then for any $\alpha>0$, $\ff$ is $\tfrac{M}{\alpha}$-relative continuous with respect to $\alpha^2\hh$.\
\item[2.] If $\ff$ is $M$-relative continuous with respect to $\hh$, then for any $\alpha>0$, $\alpha\ff$ is $M$-relative continuous with respect to $\alpha^2\hh$.\
\item[3.] If $f_j$ is $M$-relative continuous with respect to $h_j$ for $j = 1, \ldots, n$, then $\sum_{j=1}^n f_j$ is $\sqrt{n}M$-relative continuous with respect to $\sum_{j=1}^n h_j$.
\item[4.] If $f_j$ is $M_j$-relative continuous with respect to $h_j$ for $j = 1, \ldots, n$, then for $\alpha_j > 0$ and $M > 0$ it holds that $\sum_{j=1}^n \alpha_j f_j$ is $\sqrt{n}M$-relative continuous with respect to $\sum_{j=1}^n \tfrac{\alpha_j^2}{\beta_j^2} h_j$ with $\beta_j := \tfrac{M}{M_j}$.
\end{description}
\end{prop}
{\bf Proof: }   Let $x,y \in Q$, $x \ne y$, and $g(x) \in \partial f(x)$.

1. It holds that
\[
\| g(x)\|_* \le \frac{ M \sqrt{2D_h(y,x)}}{\|y-x\|} = \frac{\tfrac{M}{\alpha}  \sqrt{2D_{\alpha^2 h}(y,x)}}{\|y-x\|} \ ,
\] which establishes the result.

2. Notice that $g(x)$ is a subgradient of $f(x)$ if and only if $\alpha g(x)$ is a subgradient of $\alpha f(x)$, whereby
\[
\|\alpha g(x)\|_* = \alpha\|g(x)\|_* \le \frac{\alpha M \sqrt{2D_h(y,x)}}{\|y-x\|} = \frac{M \sqrt{2D_{\alpha^2 h}(y,x)}}{\|y-x\|} ,
\]which establishes the result.

3. Any subgradient of $\sum_{j=1}^n f_j$ at $x$ can be written as $\sum_{j=1}^n g_j(x)$ where $g_j(x) \in \partial f_j(x)$ for $j=1, \ldots, n$ (see Theorem B.21 of \cite{bertsekas}).  From the triangle inequality and the relative continuity of $f_j$ we have:
\[\begin{array}{rcl}
\left\|\displaystyle\sum_{j=1}^n g_j(x) \right\|_* \le \displaystyle\sum_{j=1}^n \left\| g_j(x) \right\|_* & \le & \displaystyle\frac{M\left(\sum_{j=1}^n \sqrt{2D_{h_j}(y,x)} \right)}{\|y-x\|}  \\ \\
&\le& \displaystyle\frac{\sqrt{n}M\left(\sqrt{2 \sum_{j=1}^n D_{h_j}(y,x) }\right)}{\|y-x\|} = \frac{\sqrt{n}M\sqrt{2 D_{h_1 + \cdots + h_n}(y,x) }}{\|y-x\|} \ ,
\end{array}\]

where the third inequality is an application of the $\ell_1$/$\ell_2$-norm inequality applied to the $n$-tuple $( \sqrt{2D_{h_1}(y,x)}, \ldots, \sqrt{2D_{h_n}(y,x)})$.

4. It follows from part (2.) that $\alpha_j f_j$ is $M_j$-continuous relative to $\alpha_j^2 h_j$. Thus $\alpha_j f_j$ is also $M$-continuous relative to $\tfrac{\alpha_j^2}{\beta_j^2} h_j$ from part (1.), whereby the proof is finished by utilizing part (3.).
\qed\medskip\medskip

We also make use of the notion of ``relative strong convexity'' which was introduced in \cite{lu2018relatively}, and is used here in some of the convergence guarantee analyses.\medskip\medskip

\begin{mydef}\label{strong}
$\ff$ is $\mu$-strongly convex relative to $\hh$ on $Q$ if there is a scalar $\mu \ge 0$ such that for any
$x,y\in \intQ$ and any $g(x) \in \partial \ff$ it holds that
\begin{equation}\label{eq:strong}
f(y) \ge f(x) + \langle g(x), y-x \rangle +  \mu \Dh(y,x) \ .
\end{equation}
\end{mydef}

In \cite{lu2018relatively} it was shown that the notion of relative strong convexity embodied in Definition \ref{strong} is the natural way to define strong convexity in the context of mirror descent and similar algorithms, and leads to linear convergence of mirror descent in the smooth setting.  In the non-smooth setting, we will show in Theorem \ref{thm:deterstrong1} that relative strong convexity improves the convergence of mirror descent from $O(1/\sqrt{k})$ to $O(1/k)$.

\subsection{Stochastic Setting}

For $x \in Q$, let $\tilde g(x)$ denote a random stochastic estimate of a subgradient of $\ff$ at $x$. Extending the definition of relative continuity from the deterministic setting, we define stochastic relative continuity as follows.\medskip

\begin{mydef}\label{cliffhanger} $\ff$ is $G$-stochastically-relative continuous with respect to the reference function  $\hh$ on $Q$ for some $G >0$ if $\ff$ together with the oracle to compute a stochastic subgradient satisfies:
\begin{enumerate}
\item Unbiasedness property: \ \  $\EE[\tg(x)|x] \in \partial f(x)$, and
\item Boundedness property: \ \  $ \EE [ \|\tg(x)\|_*^2 | x] \le G^2 \frac{D_h(y,x)}{\frac{1}{2}\|y-x\|^2} \ \ \text{for\ all} \ x, y\in Q $ and $x \ne y$.
\end{enumerate}
\end{mydef}
(In the particular case when $h(x) = \tfrac{1}{2} \|x\|_2^2$, the Bregman distance is $D_h(y,x) = \tfrac{1}{2} \|y-x\|_2^2$, whereby the stochastically-relative continuity boundedness property becomes $ \EE[\|\tg(x)\|^2_2 | x]\le G^2$ for all $x \in Q$, which corresponds to the standard condition \eqref{def:tra_stoc_cont} for the $\ell_2$-norm.) \medskip

For $x \in Q$, define
\begin{equation}
\tM(x) := \|\tg(x)\|_*\max_{ y\in Q, y\ne x } \frac{\|y-x\|}{\sqrt{2D_h(y,x)}} \ .
\end{equation}
Notice for a given $x$ that $\max_{y\in Q, y\ne x} \frac{\|y-x\|}{\sqrt{2D_h(y,x)}}$ is a deterministic quantity, and therefore $\tM(x)$ is a conditional random variable (given $x$) that is defined on the same probability space as $\tg(x)$. Clearly, if $\ff$ is $G$-stochastically-relative continuous, we have by the boundedness property that for any $x\in Q$
\begin{equation}\label{eq:boundontM}
\EE[\tM(x)^2|x]=\EE[\|\tg(x)\|_*^2|x]\max_{ y\in Q, y\ne x } \frac{\|y-x\|^2}{2D_h(y,x)}\le G^2 \ .
\end{equation}

Exactly as in the deterministic setting, we have:\medskip

\begin{prop}\label{prop:stocequiv}
If $\ff$ is $G$-stochastically-relative continuous with respect to the reference function  $\hh$ on $Q$, then for any $t>0$ it holds for all $x,y \in Q$ and any stochastic subgradient estimate $\tilde g(x)$ that:
\begin{equation*}
\tfrac{1}{t}D_h(y,x) + \langle \tg(x), y-x \rangle + \tfrac{1}{2} t\tM^2(x) \ge 0 \ .
\end{equation*}
\end{prop}
{\bf Proof:} For any $t>0$, we have
\begin{equation*}
- \langle \tg(x), y-x \rangle\le \|\tg(x)\|_* \|y-x\| \le \tM(x) \sqrt{2D_h(y,x)} \le \tfrac{1}{2} t\tM(x)^2 + \frac{D_h(y,x)}{t} \ ,
\end{equation*}
and the proof is finished by rearranging terms. \qed

\section{Computational Analysis for Stochastic Mirror Descent and (Deterministic) Mirror Descent}\label{algorithms}

In this section we present computational guarantees for Stochastic Mirror Descent (Algorithm \ref{al:proxstocgradscheme}) for minimizing a convex function $\ff$ that is $G$-stochastically-relative continuous with respect to a given reference function $\hh$.  We also present computational guarantees for (deterministic) Mirror Descent (Algorithm \ref{al:proxgradscheme}) when $\ff$ is $M$-relative continuous with respect to a reference function $\hh$, which follows as a special case of the stochastic setting.

We begin by recalling the standard Three-Point Property for optimization using Bregman distances: \medskip

\begin{lem} {\bf (Three-Point Property (Tseng \cite{tseng}))} Let $\phi(x)$ be a convex function, and let $\Dh(\cdot, \cdot)$ be the Bregman distance for $\hh$. For a given vector $z$, let
\begin{equation*}
z^+ := \arg\min_{x\in Q} \left\{ \phi(x) + \Dh(x,z) \right\} \ .
\end{equation*}
Then
\begin{equation*}
\phi(x) + \Dh(x,z) \ge \phi(z^+) + \Dh(z^+, z) + \Dh(x,z^+)\ \   for\ all \ x\in Q \ . \qed
\end{equation*}
\end{lem}\medskip

Let us denote the (primitive) random variable at the $i^{\mathrm{th}}$ iteration of the Stochastic Mirror Descent Algorithm (Algorithm \ref{al:proxstocgradscheme}) by $\gamma_{i}$, i.e., $\gamma_i$ is the random variable that determines the (stochastic) subgradient $\tilde g(x^i)$ at iterate $x^i$ in the Stochastic Mirror Descent Algorithm.  Then $x^{i+1}$ is computed according to the update of the Stochastic Mirror Descent Algorithm, whereby $x^{i+1}$ is a random variable which depends on all previous values $\gamma_0, \ldots, \gamma_i$ and we denote this string of random variables by
\[
\xi_{i} := \{\gamma_0, \ldots, \gamma_{i}\} .
\]
The following theorem states convergence guarantees for the Stochastic Mirror Descent Algorithm in terms of expectation.\medskip

\begin{thm}{\bf{(Convergence Bound for Stochastic Mirror Descent Algorithm})}\label{thm:stonon}
Consider the Stochastic Mirror Descent Algorithm (Algorithm \ref{al:proxstocgradscheme}) with given step-size sequence $\{t_i\}$. If $\ff$ is $G$-stochastically-relative continuous with respect to $\hh$ for some $G>0$, then the following inequality holds for all $k \ge 1$ and $x\in Q$:

\begin{equation}\label{stojoint}
\ \ \ \ \ \ \ \ \EE_{\xi_{k-1}}\left[ f(\barx^k)\right]  - f(x)   \ \ \leq \ \
\frac{ \tfrac{1}{2}G^{2}\sum_{i=0}^{k}t_{i}^{2}+D_h(x,x^{0})}{\sum_{i=0}^{k}t_{i}} \ ,
\end{equation}
where $\barx^k:=\frac{1}{\sum_{i=0}^k t_i} \sum_{i=0}^kt_i x^i$.
\end{thm}\medskip\medskip

{\bf Proof:}
First notice that
\begin{equation}\label{eq:stocproofmid}
\begin{array}{cl}
  f(x^{i})+ \left\langle g(x^i),x-x^{i}\right\rangle
  & = f(x^i) + \left\langle \EE_{\gamma_{i}} [\tg(x^i)| x^i],x-x^{i}\right\rangle \\ \\
 & = f(x^i) + \EE_{\gamma_{i}} \left[\left\langle  \tg(x^i),x-x^{i}\right\rangle \right | x^i] \\ \\
& \ge  f(x^{i})+\EE_{\gamma_{i}} \left[ \left\langle \tg(x^i),x^{i+1}-x^{i}\right\rangle +\frac{1}{t_{i}}D_h(x^{i+1},x^{i})+\frac{1}{t_{i}}D_h(x,x^{i+1})-\frac{1}{t_{i}}D_h(x,x^{i})  | x^i \right]\\ \\
& \ge  f(x^{i}) +  \EE_{\gamma_{i}} \left[ -\tfrac{1}{2}\tM(x^i)^{2}t_{i}+\frac{1}{t_{i}}D_h(x,x^{i+1})-\frac{1}{t_{i}}D_h(x,x^{i}) | x^i\right] \\ \\
 & \ge f(x^i) - \tfrac{1}{2} G^2 t_i + \tfrac{1}{t_i} \EE_{\gamma_{i}} [ D_h(x, x^{i+1}) | x^i] -\frac{1}{t_{i}}D_h(x,x^{i}) \ ,
\end{array}
\end{equation}

where the first equality uses the unbiasedness of $\tg(x)$, the second equality is because of linearity, the first inequality is from the Three-Point Property with $\phi(x)=t_i \langle \tg(x^i), x - x^i \rangle$, the second inequality uses Proposition \ref{prop:stocequiv}, and the last inequality uses \eqref{eq:boundontM}.  Since also $f(x) \ge  f(x^{i})+ \left\langle g(x^i),x-x^{i}\right\rangle $ from the definition of a subgradient, we have from \eqref{eq:stocproofmid}:
\[
f(x) \ge  f(x^{i})+ \left\langle g(x^i),x-x^{i}\right\rangle \ge f(x^i) - \tfrac{1}{2} G^2 t_i + \tfrac{1}{t_i} \EE_{\gamma_{i}} [ D_h(x, x^{i+1}) | x^i] -\tfrac{1}{t_{i}}D_h(x,x^{i}) \ .
\]

Taking expectation with respect to $\xi_{i}$ on both sides of the above inequality yields:
\begin{equation}
f(x) \ge \EE_{\xi_{i-1}} [f(x^i)] - \tfrac{1}{2} G^2 t_i + \tfrac{1}{t_i} \EE_{\xi_i} [ D_h(x, x^{i+1}) ] -\tfrac{1}{t_{i}}\EE_{\xi_{i-1}}[D_h(x,x^{i})]
\end{equation}
Now rearrange and multiply through by $t_{i}$ to yield:
\[
t_{i}\EE_{\xi_{i-1}}[f(x^{i}) - f(x)]\le \tfrac{1}{2}G^{2}t_{i}^{2}+ \EE_{\xi_{i-1}} [D_h(x,x^{i})]-\EE_{\xi_{i}} [D_h(x,x^{i+1})].
\]
Summing up the above inequality over $i$ and noting that $D_h(x,x^{k+1})\ge0$ we arrive at:
\begin{equation}\label{eq:lastinproof}
\begin{array}{rcl}
\tfrac{1}{2}G^{2}\sum_{i=0}^{k}t_{i}^{2}+D_h(x,x^{0}) \ge \sum_{i=0}^{k}t_{i}\EE_{\xi_{i-1}}[f(x^{i}) - f(x)] & = &  \EE_{\xi_{k-1}}\sum_{i=0}^{k}t_{i}[f(x^{i}) - f(x)] \\ \\ & \ge & \left(\sum_{i=0}^{k}t_{i}\right)\EE_{\xi_{k-1}}\left[f(\barx^{k}) - f(x) \right] \ ,
\end{array}
\end{equation}
where the last inequality uses the convexity of $\ff$.
Dividing by $\sum_{i=0}^{k}t_{i}$ completes the proof. \qed \medskip

\begin{rem}
As a direct consequence of \eqref{eq:lastinproof}, we obtain the following result which is similar to the deterministic setting \eqref{sd_bound3}:
$$
\EE_{\xi_{k-1}}\left[\min_{0\le i \le k} f(x^i)\right]   - f(x) \ \ \leq \ \
\frac{ \tfrac{1}{2}G^{2}\sum_{i=0}^{k}t_{i}^{2}+D_h(x,x^{0})}{\sum_{i=0}^{k}t_{i}} \ . $$
\end{rem}

Theorem \ref{thm:stonon} implies the following high-probability result using a simple Markov bound.\medskip

\begin{cor}
Let $x^*$ be an optimal solution of \eqref{poi1}.  Under the conditions of Theorem \ref{thm:stonon}, for any $\delta > 0$ it holds that:
\begin{equation*}
\PP\left[f(\barx^k)  - f^* \ge \delta  \right] \ \ \leq \ \
\frac{ \tfrac{1}{2}G^{2}\sum_{i=0}^{k}t_{i}^{2}+D_h(x^*,x^{0})}{\delta \sum_{i=0}^{k}t_{i}}
\end{equation*}
\end{cor}
{\bf Proof:} Using the Markov inequality, we have:
\[
\PP\left[ f(\barx^k)  - f^* \ge \delta  \right] \le \frac{\EE\left[ f(\barx^k)  - f^*  \right]}{\delta} \le \frac{ \tfrac{1}{2}G^{2}\sum_{i=0}^{k}t_{i}^{2}+D_h(x^*,x^{0})}{\delta \sum_{i=0}^{k}t_{i}} \ .
\]
\qed \medskip

Similar to the case of traditional analysis of stochastic mirror descent, the Stochastic Mirror Descent Algorithm (Algorithm \ref{al:proxstocgradscheme}) leads to an $O(\tfrac{1}{\varepsilon^2})$ convergence guarantee (in expectation) by using an appropriate step-size sequence $\{t_i\}$ as the next corollary shows. \medskip

\begin{cor}\label{cor:convergence1}
Under the conditions of Theorem \ref{thm:stonon}, for a given $\varepsilon > 0$ suppose that the step-sizes are set to:
\[
t_i:= \frac{\varepsilon}{G^2} \
\]
for all $i$.  Then within
\[
k:= \left\lceil \frac{2G^2 D_h(x^*, x^0)}{\varepsilon^2} \right\rceil - 1
\]
iterations of the Stochastic Mirror Descent Algorithm it holds that:
\[
\EE\left[f(\barx^k)  \right]  - f^* \le \varepsilon \ ,
\]
where $x^*$ is any optimal solution of \eqref{poi1}.

\end{cor}
{\bf Proof:} Substituting the values of $t_i$ in \eqref{stojoint} yields the result directly. \qed \medskip

\begin{rem} Similar to the standard stochastic gradient descent scheme, the step-size rule in Corollary \ref{cor:convergence1} leads to the optimal rate of convergence provided by the bound in Theorem \ref{thm:stonon}.
\end{rem}\medskip

\begin{rem}\label{rem:nonstrong_comparison} Let us now compare these results to related results of Nedi\'{c} and Lee \cite{nedic2014stochastic}.  In order to attain an $\varepsilon$-optimality gap, \cite{nedic2014stochastic} proved a bound of
$
\left\lceil\frac{9G_f^2D_{\max}}{4\mu_h\varepsilon^2}\right\rceil
$
iterations, which follows by rearranging \eqref{sto_st_nonstrong_bound}.  In addition to not requiring Lipschitz continuity of $\ff$, our bound does not require that $\hh$ be strongly convex.  We also do not require boundedness of the feasible region; and in most settings $D_h(x^*, x^0) \ll D_{\max}$ even when $D_{\max} < +\infty$.  Furthermore, even in the setting of \eqref{sto_st_nonstrong_bound}, it holds that:
\[
G^2=\EE\left[\tg(x)^2|x\right]\max_{y\in Q, y\neq x}\frac{\|y-x\|^2}{2D_h(y,x)}\le \EE\left[\|\tg(x)\|_*^2|x\right]\tfrac{1}{\mu_h}\le \tfrac{G_{f}^2}{\mu_h} \ ,
\]
where the first inequality utilizes the strong convexity (in the standard sense) of $\hh$, and the second inequality is due to the assumption that $\ff$ is $G_f$-stochasticlly continuous (in the standard sense). Thus we see that the bound in Corollary \ref{cor:convergence1} improves on the bound in \cite{nedic2014stochastic}.
\end{rem}

In the case when $\ff$ is also $\mu_f$-strongly convex relative to $\hh$ (see Definition \ref{strong}), we obtain an $O(\tfrac{1}{k})$ convergence guarantee in expectation, which is also similar to the traditional case of stochastic gradient descent. This is shown in the next result.\medskip

\begin{thm}{\bf{(Convergence Bound for Stochastic Mirror Descent Algorithm under Strong Convexity relative to $\hh$})}\label{thm:deterstrong1} Consider the Stochastic Mirror Descent Algorithm (Algorithm \ref{al:proxstocgradscheme}). If $\ff$ is $G$-stochastically-relative continuous with respect to $\hh$ for some $G>0$ and $\ff$ is $\mu$-strongly convex relative to $\hh$ for some $\mu>0$, and if the step-sizes are chosen as $t_i = \tfrac{2}{\mu (i+1)}$, then the following inequality holds for all $k \ge 1$:
\begin{equation*}
\EE_{\xi_{k-1}}\left[f(\hx^{k})\right]-f^{*}\le\frac{2G^{2}}{\mu(k+1)} \ ,
\end{equation*}
where $\hx^k:=\frac{2}{k(k+1)} \sum_{i=0}^k i \cdot x^i$.
\end{thm}

{\bf Proof:}
For any $x \in Q$ it follows from the definition of $\mu$-strong convexity \eqref{eq:strong} that
\[
f(x)  \ge  f(x^{i})+\left\langle g(x^i),x-x^{i}\right\rangle +\mu D_h(x,x^{i}) \ .
\]
Combining the above inequality with \eqref{eq:stocproofmid} yields
\begin{equation*}
f(x)  \ge  f(x^{i})-\tfrac{1}{2}G^{2}t_{i}+\tfrac{1}{t_{i}}\EE_{\gamma_{i}} [D_h(x,x^{i+1})|x^i]-(\tfrac{1}{t_{i}}-\mu)D_h(x,x^{i}) \ .
\end{equation*}
Substituting $t_{i}=\tfrac{2}{\mu(i+1)}$ and multiplying by $i$ in the above inequality yields:
\begin{equation*}
\begin{array}{lcl}
i\left(f(x^{i})-f(x)\right) & \le & \frac{G^{2}i}{\mu(i+1)}+\frac{\mu}{2}\left(i(i-1)D_h(x,x^{i})-i(i+1)\EE_{\gamma_{i}} [D_h(x,x^{i+1})|x^i]\right) \\ \\
& \le & \frac{G^{2}}{\mu}+\frac{\mu}{2}\left(i(i-1)D_h(x,x^{i})-i(i+1)\EE_{\gamma_{i}} [D_h(x,x^{i+1})|x^i]\right) \ .
\end{array}
\end{equation*}

Taking expectation over $\xi_{i-1}$ and summing up the above inequality over $i$ then yields
\[
\left(\sum_{i=1}^k i \right)\EE_{\xi_{k-1}} [f(\hx^k)-f(x)]\le \sum_{i=1}^{k}i\left(\EE_{\xi_{i-1}}[f(x^{i})]-f(x)\right) \le \frac{kG^{2}}{\mu} - k(k+1)\left(\frac{\mu}{2}\right)\EE_{\xi_{k}} [D_h(x,x^{k+1})]  \le  \frac{kG^{2}}{\mu} \ ,
\]
where the first inequality uses the convexity of $\ff$ and the observation that $\EE_{\xi_{i-1}} f(x^i) = \EE_{\xi_{k-1}} f(x^i)$ for $i \le k$.
Taking $x = x^*$ where $x^*$ is an optimal solution of \eqref{poi1}, the proof is completed by noticing $\sum_{i=1}^k i = \tfrac{k(k+1)}{2}$. \qed \medskip\medskip

\begin{rem}
It may not be easy to find cases when the objective function is both $G$-stochastically-relative continuous and is $\mu$-relatively strongly convex relative to $h$.  However, as long as these properties are satisfied along the path of iterates or around the minimum, one can achieve the faster convergence of Theorem \ref{thm:deterstrong1}.
\end{rem}\medskip

\begin{rem}
Let us also compare the computational guarantee of Theorem \ref{thm:deterstrong1} to the results in Nedi\'{c} and Lee \cite{nedic2014stochastic}.  In order to attain an $\varepsilon$-optimality gap, \cite{nedic2014stochastic} proved the bound \eqref{sto_st_strong_bound}.  First notice that we do not require either that $\ff$ is uniformly Lipschitz continuous or that $\hh$ is strongly convex in the traditional sense, or that $h$ is uniformly smooth.  However, even if these requirements hold, it follows from Remark \ref{rem:nonstrong_comparison} that $G^2\le \frac{G_f^2}{\mu_h}$, and it also holds that:
\[
D_f(x,y)\ge \tfrac{\mu_f}{2}\|x-y\|^2\ge\tfrac{\mu_f}{L_h}D_h(x,y) \ ,
\]
where the first inequality utilizes that $\ff$ is $\mu_f$ strongly convex and the second inequality $\hh$ is $L_h$ smooth in the standard sense. Thus $\ff$ is at least $\mu=\tfrac{\mu_f}{L_h}$-strongly convex relative to $\hh$ (this follows by applying Proposition 1.1 in \cite{lu2018relatively}). Therefore, even under the stronger requirements of \cite{nedic2014stochastic}, Theorem \ref{thm:deterstrong1} improves on the corresponding result in \cite{nedic2014stochastic}.

\end{rem}
We end this section with a discussion of the deterministic setting, namely the (Deterministic) Mirror Descent Algorithm (Algorithm \ref{al:proxgradscheme}).  Suppose that there is no stochasticity in the computation of subgradients.  We can cast this as an instance of the Stochastic Mirror Descent Algorithm (Algorithm \ref{al:proxstocgradscheme}) wherein $\tg(x) = g(x) \in \partial f(x)$ for all $x \in Q$.  In this case relative stochastic continuity (Definition \ref{cliffhanger})  is equivalent to relative continuity (Definition \ref{def:continuity}) with the same constant. Thus deterministic Mirror Descent is a special case of Stochastic Mirror Descent, and we have the following computational guarantees as special cases of the stochastic case.\medskip

\begin{thm}{\bf{(Convergence Bound for Deterministic Mirror Descent Algorithm})}\label{thm:deternon}
Consider the (Deterministic) Mirror Descent Algorithm (Algorithm \ref{al:proxgradscheme}). If $\ff$ is $M$-relative continuous with respect to $\hh$ for some $M>0$, then for all $k \ge 1$ and $x\in Q$ the following inequality holds:
\begin{equation*}\label{joint}
\ \ \ \ \ \ \  f(\barx^k)  - f(x)   \ \ \leq \ \
 \frac{\frac{1}{2}M^{2}\sum_{i=0}^{k}t_{i}^{2}+D_h(x,x^{0})}{\sum_{i=0}^{k}t_{i}} \ ,
\end{equation*}
where $\barx^k:=\frac{1}{\sum_{i=0}^k t_i} \sum_{i=0}^kt_i x^i$.
\end{thm} \qed \medskip

\begin{cor}\label{cor:convergence2}
Under the conditions of Theorem \ref{thm:deternon}, for a given $\varepsilon > 0$ suppose that the step-sizes are set to:
\[
t_i:= \frac{\varepsilon}{M^2} \
\]
for all $i$.  Then within
\[
k:= \left\lceil \frac{2M^2 D_h(x^*, x^0)}{\varepsilon^2} \right\rceil - 1
\]
iterations of Deterministic Mirror Descent it holds that:
\[
f(\barx^k)  - f^* \le \varepsilon \ ,
\]
where $\barx^k:=\frac{1}{\sum_{i=0}^k t_i} \sum_{i=0}^kt_i x^i$,  and $x^*$ is any optimal solution of \eqref{poi1}.
\end{cor} \qed \medskip

\begin{thm}{\bf{(Convergence Bounds for Deterministic Mirror Descent with Strong Relative Convexity})}\label{thm:deterstrong2}
Consider the Deterministic Mirror Descent Algorithm (Algorithm \ref{al:proxgradscheme}). If $\ff$ is $M$-relative continuous with respect to $\hh$ for some $M>0$ and $\ff$ is $\mu$-strongly convex relative to $\hh$ for some $\mu>0$, and if the step-sizes are chosen as $t_i = \tfrac{2}{\mu (i+1)}$, then the following inequality holds for all $k \ge 1$:
\begin{equation*}
f(\hx^{k})-f^{*}\le\frac{2M^{2}}{\mu(k+1)} \ ,
\end{equation*}
where $\hx^k:=\frac{2}{k(k+1)} \sum_{i=0}^k i \cdot x^i$.

\end{thm} \qed

\section{Specifying a Reference Function $\hh$ with Relative Continuity for Mirror Descent}\label{wellfleet}

Let us discuss using either deterministic or stochastic Mirror Descent (Algorithm \ref{al:proxgradscheme} or Algorithm \ref{al:proxstocgradscheme}) for solving the optimization problem \eqref{poi1} with objective function $\ff$ that is  $M$-relative continuous or $G$-stochastically-relative continuous (respectively) with respect to the reference function $\hh$.  In order to efficiently execute the update step in Algorithm \ref{al:proxgradscheme} and/or Algorithm \ref{al:proxstocgradscheme} we need $\hh$ to be such that the linearization subproblem LS \eqref{subproblem} is efficiently solvable for any given $c$.  Therefore, in order execute Mirror Descent for solving \eqref{poi1} using Algorithm \ref{al:proxgradscheme} or Algorithm \ref{al:proxstocgradscheme}, we need to specify a differentiable convex reference function $\hh$ that has the following two properties:
\begin{enumerate}
\item[(i)] $\ff$ is $M$-relative continuous (or $G$-stochastically-relative continuous) with respect to $\hh$ on $Q$ for $M$ (or $G$) that is easy to determine, and
\item[(ii)] the linearization subproblem LS \eqref{subproblem} has a solution, and the solution is efficiently computable.
\end{enumerate}

We now discuss quite broadly how to construct such a reference function $\hh$ with these two properties when $\|g(x)\|_*^2$ is bounded by a polynomial in $\|x\|_2$.

\subsection{Deterministic Setting}\label{ferry}

Suppose that $\| g(x) \|_*^2 \le p_r(\|x\|_2)$ for all $x \in Q$ and all $g(x) \in \partial f(x)$, where $p_r(\alpha)=\sum_{i=0}^r a_i \alpha^i$ is an $r$-degree polynomial of $\alpha$  whose coefficients $\{a_i\}$ are nonnegative.  Let $$h(x) := \sum_{i=0}^r \tfrac{a_i}{i+2} \|x\|_2^{i+2}  \ . $$   Then the following proposition states that $\ff$ is $1$-relative continuous with respect to $\hh$.  This implies that no matter how fast the subgradient of $\ff$ grows polynomially as $\|x\|_2 \rightarrow \infty$, $\ff$ is relatively continuous with respect to the simple reference function $\hh$, even though $f$ does not exhibit uniform Lipschitz continuity. \medskip\medskip

\begin{prop}\label{lem:continuity-Rn}
$\ff$ is $1$-continuous relative to $h(x) = \sum_{i=0}^r \frac{a_i}{i+2} \|x\|_2^{i+2}$.
\end{prop}

{\bf Proof:}
Let $h_i(x) = \tfrac{1}{i+2} \|x\|_2^{i+2}$, then $h(x) = \sum_{i=0}^r a_i h_i(x)$, and by the definition of Bregman distance, we have
\begin{equation*}
\begin{array}{lcl}
D_{h_i}(y,x) & = &\frac{1}{i+2}\|y\|_{2}^{i+2}-\frac{1}{i+2}\|x\|_{2}^{i+2}-\left\langle \|x\|_{2}^{i} x,y-x\right\rangle \\ \\
 & = &\frac{1}{i+2}\left(\|y\|_{2}^{i+2}+(i+1)\|x\|_{2}^{i+2}-(i+2)\|x\|_{2}^{i}\left\langle x,y\right\rangle \right) \ .
\end{array}
\end{equation*}
Notice that
\begin{equation*}
\begin{array}{lcl}
 && \|y\|_{2}^{i+2}+(i+1)\|x\|_{2}^{i+2}-(i+2)\|x\|_{2}^{i}\left\langle x,y\right\rangle \\ \\
&= & \left(\|y\|_{2}^{i+2}+\frac{i}{2} \|x\|_{2}^{i+2}-\frac{i+2}{2}\|x\|_{2}^{i}\|y\|_{2}^{2}\right)+\frac{i+2}{2}\|x\|_{2}^{i}\left(\|x\|_2^2+\|y\|_2^2-2\left\langle x,y\right\rangle \right)\\ \\
& \ge & 0 + \frac{i+2}{2}\|x\|_{2}^{i}\left(\|x\|_2^2+\|y\|_2^2-2\left\langle x,y\right\rangle \right) \\ \\
& = & \frac{i+2}{2}\|x\|_{2}^{i}\|y-x\|_{2}^{2} \ ,
\end{array}
\end{equation*}
where the inequality above is an application of arithmetic-geometric mean inequality $a^\lambda b^{1-\lambda} \le \lambda a + (1-\lambda) b$ with $a = \|x\|_2^{i+2}$, $b = \|y\|_2^{i+2}$, and $\lambda = \tfrac{i}{i+2}$.  Thus we have
\begin{equation}\label{eq:boundD}
D_{h}(y,x)= \sum_{i=0}^r a_i D_{h_i}(y,x) \ge \tfrac{1}{2}\|y-x\|_{2}^{2}\left(\sum_{i=0}^r a_i \|x\|_2^{i} \right) \ .
\end{equation}

Therefore
\[
\|g(x)\|_*^2 \le p_r(\|x\|_2) = \sum_{i=0}^r a_i \|x\|_2^{i} \le \frac{D_h(y,x)}{\tfrac{1}{2}\|y-x\|_2^2} \ ,
\]
which shows that $\ff$ is $1$-relative continuous with respect to $\hh$.\qed \medskip\medskip

\noindent {\bf Solving the linearization subproblem \eqref{subproblem}.}  Let us see how we can solve the linearization subproblem \eqref{subproblem} for this class of optimization problems.  The linearization subproblem \eqref{subproblem} can be written as
\begin{equation}\label{eq:rnsubprob}
\mathrm{LS:} \ \ \ \min_{x \in \mathbb{R}^n} \ \  \langle c, x \rangle + \sum_{i=0}^r \tfrac{a_i}{i+2} \|x\|_2^{i+2} \ ,
\end{equation}
and the first-order optimality condition is simply:
\begin{equation}\label{eq:rnsubprob33}
c+ \left( \sum_{i=0}^r a_i \|x\|_2^{i} \right) x = 0 \ ,
\end{equation}
whereby $x=-\theta c$ for some scalar $\theta \ge 0$, and it remains to simply determine the value of the nonnegative scalar $\theta$.  In the case when $c=0$ we have $x=0$ satisfies \eqref{eq:rnsubprob33}, so let us examine the case when $c \ne 0$, in which case from \eqref{eq:rnsubprob33} $\theta$ must satisfy:
$$
\sum_{i=0}^r a_i \|c\|_2^{i} \theta ^{i + 1} - 1 = 0 \ ,
$$
which implies that $\theta$ is the unique positive root of a univariate polynomial monotone in $\theta \ge 0$.  For $r \in \{0, 1, 2, 3\}$ this root can be computed in closed form.  Otherwise the root can be computed efficiently (up to machine precision) using any suitable root-finding method. \medskip

\begin{rem} We can incorporate a simple set constraint $x\in Q$ in problem \eqref{eq:rnsubprob} provided that we can easily compute the Euclidean projection on $Q$. In this case, the linearization subproblem \eqref{subproblem} can be converted to a $1$-dimensional convex optimization problem, see Appendix A.1 of \cite{lu2018relatively} for details.
\end{rem}

\subsection{Stochastic Setting}\label{sec:Rnstochastic}
In the stochastic setting, the stochastic subgradient $\tg(x)$ is a conditional random variable for a given $x$. Suppose that $\EE\left[\|\tg(x)\|^2_*|x\right] \le p_r(\|x\|_2)$ for all $x \in Q$, where $p_r(\alpha) = \sum_{i=0}^r a_i \alpha^i$ is an $r$-degree polynomial whose coefficients $\{a_i\}$ are nonnegative.  Let
\[
h(x):=\sum_{i=0}^r \tfrac{a_i}{i+2}\|x\|_2^{i+2} \ ,
\] and similar to the deterministic case we have:\medskip

\begin{prop}\label{lem:continuity-Rn-friday}
$\ff$ is $1$-stochastically continuous relative to $h(x) = \sum_{i=0}^r \frac{a_i}{i+2} \|x\|_2^{i+2}$.
\end{prop}

{\bf{Proof:}} For any $x,y\in Q$ with $x \ne y$ we have:
\begin{equation}\label{eq:boundstoM}
\EE\left\|\tg(x)\|^2_* | x\right]\le p_r(\|x\|_2) =  \sum_{i=0}^r a_i \|x\|_2^i \le \frac{2D_h(y,x)}{\|y-x\|_2^2} \ ,
\end{equation}
where the last inequality follows from \eqref{eq:boundD}, and thus $\ff$ is $1$-stochastically continuous relative to $\hh$. \qed\medskip

{\bf Solving the linearization subproblem \eqref{subproblem}}.  The linear optimization subproblem is identical in structure to that in the deterministic case and so can be solved as discussed at the end of Section \ref{ferry}.

%
%

\subsection{Relative Continuity for instances of SVM and IEP}\label{wrightwoods}
Here we examine in detail the two motivating examples stated in the Introduction, namely the Support Vector Machine (SVM) problem, and the Intersection of Ellipsoids Problem (IEP).  We first prove the following lemma, which presents upper bounds on the Bregman distances $D_h(
\cdot, \cdot)$ for $h(x)=\frac{1}{3}\|x\|_2^3$ and $h(x)=\frac{1}{4}\|x\|_2^4$.\medskip

\begin{lem}\label{thm:lem_for_4th}
$ $
\begin{enumerate}

\item Let $h(x):=\frac{1}{3}\|x\|_2^3$.  Then $D_h(y,x)\le \frac{1}{3}\|y-x\|^2_2\left(\|y\|_2+2\|x\|_2\right)$.
\item Let $h(x):=\frac{1}{4}\|x\|_2^4$.  Then $D_h(y,x)\le \frac{1}{4}\|y-x\|^2_2\left(\|y+x\|_2^2+2\|x\|_2^2\right)$.

\end{enumerate}
\end{lem}
{\bf Proof:}

1.
\begin{equation*}
\begin{array}{lcl}
D_h(y,x)&=&\frac{1}{3}\left(\|y\|_2^3+2\|x\|_2^3-3\|x\|_2\langle x,y\rangle\right)\\ \\
&\le& \frac{1}{3}\left(\|y\|_2^3+2\|x\|_2^3-3\|x\|_2\langle x,y\rangle -2\|y\|_2\langle y,x \rangle + 2\|y\|_2^2\|x\|_2- \|x\|_2\langle x, y \rangle + \|y\|_2\|x\|_2^2\right) \\ \\
&=& \frac{1}{3}\|y-x\|^2_2\left(\|y\|_2+2\|x\|_2\right) \ ,
\end{array}
\end{equation*}
where the first equality follows from simplifying and combining terms, the inequality follows from applying the Cauchy-Schwarz inequality twice, and the final equality is from simplifying and combining terms.\medskip

2.
\begin{equation*}
\begin{array}{lcl}
D_h(y,x)&=&\frac{1}{4}\left(\|y\|_2^4+3\|x\|_2^4-4\|x\|_2^2\langle x,y\rangle\right)\\ \\
&\le& \frac{1}{4}\left(\|y\|_2^4+3\|x\|_2^4-4\|x\|_2^2\langle x,y\rangle + 4\|x\|_2^2\|y\|_2^2 - 4\langle x,y\rangle^2\right) \\ \\
&=& \frac{1}{4}\|y-x\|^2_2\left(\|y+x\|_2^2+2\|x\|_2^2\right) \ ,
\end{array}
\end{equation*}
where the first equality follows from simplifying and combining terms, the inequality follows from applying the Cauchy-Schwarz inequality once, and the final equality is from simplifying and combining terms.\qed \medskip

{\bf Support Vector Machine (SVM).} The Support Vector Machine (SVM) is an important supervised learning model for binary classification in machine learning.  The SVM optimization problem for binary classification is:
\begin{equation}\label{eq:svm}
\mbox {SVM:} \ \ \ \min_x f(x):= \tfrac{1}{n} \sum_{i=1}^n \max\{0, 1-y_i x^T w_i\} + \tfrac{\lambda}{2}\|x\|_2^2 \ ,
\end{equation}

where $w_{i}$ is the input feature vector of sample $i$ and $y_{i} \in \{-1,1\}$ is the label of sample $i$.  Notice that $\ff$ is not differentiable due to the presence of the hinge loss terms in the summation, and $\ff$ is also not Lipschitz continuous due to the presence of the $\ell_2$-norm regularization term; thus we cannot directly utilize typical subgradient or gradient schemes and their associated computational guarantees in the analysis of \eqref{eq:svm}.  Researchers have developed various approaches to overcome this limitation.  For example, \cite{duchi2009efficient} introduced a splitting subgradient-type method, where the basic idea is to split the loss function and the regularization terms.  \cite{yu2010quasi} introduced a quasi-Newton method, where they do not need to worry about the unbounded subgradient. Another approach is to {\em a priori} constrain $x$ to lie in an $\ell_2$-ball of radius $R$ for $R$ sufficiently large so that the ball contains the optimal solution, and to project onto this ball at each iteration; in this approach $\ff$ is Lipschitz continuous in the amended feasible region, see \cite{shalev2007pegasos}.   Indeed, one can show using quadratic optimization optimality conditions that it suffices to set $R = \min\{\tfrac{1}{\lambda}\left(\tfrac{1}{n}\sum_{i=1}^n \|w_i\|_2\right), \sqrt{2/\lambda}\}$ (see Appendix \ref{sec:appendix}) wherein the modulus of Lipschitz continuity in the amended feasible region is at most $M \le \tfrac{1}{n}\sum_{i=1}^n \|w_i\|_2 + \min\{\sqrt{2\lambda}, \tfrac{1}{n}\sum_{i=1}^n \|w_i\|_2 \}$. Furthermore, in \cite{lacoste2012simpler} the authors show that if the initial point lies within a suitably chosen large ball, then Stochastic Subgradient Descent with a small step-size ensures in expectation that all iterates lie in the large ball, which then ensures that the norms of all subgradients are bounded in expectation.  

Let us see how we can directly use the constructs of relative continuity to tackle the SVM problem with a suitably designed version of Stochastic Mirror Descent -- without any projection step to a ball.  We can rewrite the objective function of \eqref{eq:svm} as $$f(x) = \tfrac{1}{n} \sum_{j=1}^n f_j(x) \ ,$$ where $f_j(x) = \max\{0, 1-y_j x^T w_j\} + \frac{\lambda}{2}\|x\|_2^2$. We consider computing a stochastic estimate of the subgradient of $\ff$ by using a single sample index $\tilde j$ drawn randomly from $\{1, \ldots, n\}$, namely $\tg(x) \in \partial f_{\tilde j}(x)$ where $\tilde j$ is drawn uniformly at random from $\{1,\ldots,n\}$.  Then $\| \tg(x) \|_2^2 \le  (\lambda \|x\|_2 + \|w_{\tilde j}\|_2)^2$, whereby
\[
\EE[\|\tg(x)\|_2^2|x]\le \lambda^2\|x\|_2^2+\frac{2\lambda}{n}\left(\sum_{i=1}^n \|w_i\|_2\right)\left\| x\right\|_2+ \frac{1}{n}\sum_{i=1}^n \|w_i\|_2^2 \ ,
\]
and notice that the right-hand side is a polynomial in $\|x\|_2$ of degree $r=2$.  If we choose the reference function $\hh$ as \begin{equation}\label{hyannis} h(x) := \frac{\lambda^2}{4} \|x\|_2^4 + \frac{2\lambda}{3n}\left(\sum_{i=1}^n \|w_i\|_2\right) \|x\|_2^3 + \frac{1}{2n}\left(\sum_{i=1}^n \|w_i\|_2^2 \right)\|x\|_2^2 \ , \end{equation} it follows from the Proposition \ref{lem:continuity-Rn-friday} that $\ff$ is $1$-stochastically continuous relative to $h(x)$.\medskip\medskip

\begin{prop} {\bf (Computational Guarantees for Stochastic Mirror Descent for the SVM problem \eqref{eq:svm}.)}\label{hampton}  Consider applying the Stochastic Mirror Descent algorithm (Algorithm \ref{al:proxstocgradscheme}) to the Support Vector Machine problem \eqref{eq:svm} using the reference function \eqref{hyannis}.  For an absolute optimality tolerance value $\varepsilon >0$, and using the constant step-sizes $t_i:=\varepsilon$, let the algorithm be run for
\small{\[ k :=
\left\lceil \frac{\|x^*-x^0\|^2 \left(3\lambda^2\left(\|x^* + x^0\|_2^2 + 2
\|x^0\|_2^2\right) + \frac{8\lambda}{n}\left(\sum_{i=1}^n\|w_i\|_2\right)\left(\|x^*\|_2 + 2\|x^0\|_2\right)+ \frac{6}{n}\left(\sum_{i=1}^n\|w_i\|_2^2\right) \right) }{6\varepsilon^2}  \right\rceil -1
\]
}\normalsize{
iterations, where $x^*$ is the optimal solution of \eqref{eq:svm}.  Then it holds that
$$\EE\left[ f(\barx^k) - f^*\right] \le \varepsilon \ , $$
where $\barx^k:=\frac{1}{k+1}\sum_{i=0}^k x^i$.}
\end{prop}

\noindent {\bf Proof:}  We showed above (using Proposition \ref{lem:continuity-Rn-friday}) that $\ff$ is $1$-stochastically continuous relative to $\hh$ defined in \eqref{hyannis}.  Furthermore, applying Lemma \ref{thm:lem_for_4th} it follows that \small{$$D_h(x^*,x^0)\le  \frac{\lambda^2}{4}\|x^* - x^0\|^2_2\left(\|x^*+x^0\|_2^2+2\|x^0\|_2^2\right)+\frac{2\lambda}{3n}\left(\sum_{i=1}^n\|w_i\|_2\right)\left(\|x^*\|_2 + 2\|x^0\|_2\right)+ \frac{1}{2n}\left(\sum_{i=1}^n\|w_i\|_2^2\right)\|x^*-x^0\|_2^2\ .$$}\normalsize{The proof is finished by substituting these values into the computational guarantee of Corollary \ref{cor:convergence1}.} \qed

{\bf Intersection of Ellipsoids Problem (IEP).\footnote{This problem was suggested by Nesterov \cite{nesterov-private}.}} Consider the problem of computing a point $x\in \RR^m$ in the intersection of $n$ ellipsoids, namely:
\begin{equation}\label{brook}
x \in \cQ := \cQ_1 \cap \cQ_2 \cap \cdots \cap \cQ_n \ ,
\end{equation}
where $\cQ_i= \{x\in \RR^m : \tfrac{1}{2}x^T A_i x + b_i x + c_i \le 0\}$ and $A_i\in \RR^{m\times m}$ is a given symmetric positive semi-definite matrix, $i = 1, \ldots, n$.  This problem can be cast as a second-order cone optimization problem, and hence can be tackled using interior-point methods.  However, interior-point methods are typically only effective when the dimensions $m$ and/or $n$ are of moderate size.  On the other hand, another way to tackle the problem is to use a first-order method to solve the unconstrained problem
\begin{equation}\label{pro:quadr}
\mbox{IEP:} \ \ \ f^* = \min_x f(x) := \max_{0 \le i \le n} \{ \tfrac{1}{2} x^T A_i x + b_i^T x + c_i \}\ ,
\end{equation}
and notice that $f(x) \le 0  \Leftrightarrow x \in \cQ$, and $\cQ \ne \emptyset \Leftrightarrow f^* \le 0$.  However, the objective function $\ff$ in \eqref{pro:quadr} is both non-differentiable and non-Lipschitz, and so it falls outside of the scope of optimization problems for which traditional first-order methods are applicable. Let us see how we can use the machinery of relative continuity to tackle this problem.  Let $\sigma := \max_{0 \le i \le n}\|A_i\|_2^2$ where $\|A_i\|_2$ is the spectral radius of $A_i$, let $\rho:=2 \max_{0 \le i \le n} \|A_i b_i\|_2$ and let $\gamma := \max_{0 \le i \le n} \|b_i\|_2^2$.  Notice that for any $x$ and $i=1, \ldots, n$, we have $g_i(x) := A_i x + b_i =\nabla f_i(x)$ where $f_i(x)$ is the $i^{\mathrm{th}}$ term in the objective function of \eqref{pro:quadr}.  Since $g(x) \in \partial f(x)$ if and only if $g(x)$ is a convex combination of the active gradients $\nabla f_i(x)$ (see Danskin's Theorem, Proposition B.22 in \cite{bertsekas}), it follows for any $g(x) \in \partial f(x)$ that
\[
\|g(x)\|_2^2 \le \max_{0 \le i \le n}\|A_i x + b_i\|^2_2 \le \max_{0 \le i \le n}\|A_i\|^2_2 \|x\|_2^2 + 2\|b_i^TA_i\|_2\|x\|_2 + \|b_i\|_2^2 \le \sigma\|x\|_2^2 + \rho\|x\|_2 + \gamma \ .
\]
Therefore we have $\|g(x)\|^2_2 \le p_2(\|x\|_2)$, where $p_1(\alpha) = \sigma \alpha^2 + \rho \alpha + \gamma$ is a quadratic function of $\alpha$, which is a polynomial in $\alpha$ of degree $r=2$.  It follows from Proposition \ref{lem:continuity-Rn} that $\ff$ is $1$-continuous relative to the reference function
\begin{equation}\label{boulevard}
h(x) := \tfrac{\sigma}{4} \|x\|_2^4 + \tfrac{\rho}{3} \|x\|_2^3+ \tfrac{\gamma}{2} \|x\|_2^2 \ .
\end{equation}

\begin{prop} {\bf (Computational Guarantees for Deterministic Mirror Descent for the IEP problem \eqref{pro:quadr}).}\label{bridge}  Consider applying the Deterministic Mirror Descent algorithm (Algorithm \ref{al:proxgradscheme}) to the Ellipsoid Intersection Problem \eqref{pro:quadr} using the reference function \eqref{boulevard}, where $\sigma := \max_{0 \le i \le n}\|A_i\|_2^2$ and $\|A_i\|_2$ is the spectral radius of $A_i$, $\rho:=2 \max_{0 \le i \le n} \|A_i b_i\|_2$ and $\gamma := \max_{0 \le i \le n} \|b_i\|_2^2$.  For an absolute optimality tolerance value $\varepsilon >0$, and using the constant step-sizes $t_i:=\varepsilon$, let the algorithm be run for
 \[
k :=\left\lceil \frac{\|x^*-x^0\|^2 \left(3\sigma\left(\|x^* + x^0\|_2^2 + 2
\|x^0\|_2^2\right) + 4\rho\left(\|x^*\|_2 + 2\|x^0\|_2\right)+ 6 \gamma\right) }{6\varepsilon^2}  \right\rceil -1
\]
 iterations, where $x^*$ is any optimal solution of \eqref{pro:quadr}.  Then it holds that
$$ f(\barx^k) - f^* \le \varepsilon \ , $$
where $\barx^k:=\frac{1}{k+1}\sum_{i=0}^k x^i$.
\end{prop}

\noindent {\bf Proof:}  We showed above (using Proposition \ref{lem:continuity-Rn}) that $\ff$ is $1$-continuous relative to $h(x)=\frac{\sigma}{4}\|x\|_2^4+\frac{\rho}{3}\|x\|_2^3+\frac{\gamma}{2}\|x\|_2^2$.  Furthermore, applying Lemma \ref{thm:lem_for_4th} it follows that $D_h(x^*,x^0)\le \frac{\sigma}{4}\|x^*-x^0\|^2_2\left(\|x^*+x_0\|_2^2+2\|x_0\|_2^2\right)+\frac{\rho}{3}\|x^*-x^0\|_2^2(\|x^*\|_2+2\|x^0\|_2)+\frac{\gamma}{2}\|x^*-x^0\|_2^2$.  The proof is finished by substituting these values into the computational guarantee of Corollary \ref{cor:convergence2}. \qed\medskip


\section*{Acknowledgement}
The author would like to express his gratitude to Robert M. Freund for thoughtful discussions that helped motivate this work, for commenting on earlier drafts of this paper, and for advising on the presentation and positioning of this paper. The author also wishes to thank Yurii Nesterov for encouraging the author's work on this topic, and for pointing out the application of IEP.

\section*{Appendix:  Finite Radius Bound for SVM}\label{sec:appendix}
Here we derive an upper bound on the norm of an optimal solution of the SVM problem \eqref{eq:svm}.\medskip

\begin{prop}
The optimal solution to the SVM problem \eqref{eq:svm} lies in the ball $B_2(0,R)$ for $R=\min\left\{\frac{1}{n\lambda}\sum_{i=1}^n\|w_i\|_2, \sqrt{2/\lambda}\right\}$.
\end{prop}

{\bf Proof:} For convenience define $A_i := y_i w_i$ for $i=1, \ldots, n$.  Then we can re-write the SVM problem as the following constrained optimization problem:

$$\begin{array}{rrl}
\min_{s,x} & \tfrac{1}{n}e^Ts + \tfrac{\lambda}{2}\|x\|_2^2 \\ \\
\mbox{s.t.} & s + Ax \ge e \\ \\
& s \ge 0 \ . \\
\end{array}$$
Let $\pi$ and $\beta$ be the multipliers on the inequality constraints above.  Then the KKT conditions imply, among other things, that the optimal solution $x^*$ must satisfy:

$$\begin{array}{rcl}
\pi^* + \beta^* = \tfrac{1}{n}e \\ \\
\lambda x^* = A^T\pi^*
\end{array}$$

where $\pi^* \ge 0$ and $\beta^* \ge 0$.  Define $\bar\pi^* =  n\pi^*$.  Then $0 \le \bar\pi^* \le e$ and

$$\lambda \|x^*\|_2 =  \|A^T\pi^*\|_2 = \tfrac{1}{n}\|A^T\bar\pi^* \|_2 \le  \tfrac{1}{n}\sum_{i=1}^n\|A_i\|_2  = \tfrac{1}{n}\sum_{i=1}^n\|w_i\|_2 \ , $$

which proves the first term in the definition of $R$.  Also, we have $
\tfrac{\lambda}{2}\|x^*\|_2^2 \le f(x^*) \le f(0) = 1,
$
thus $\|x^*\|_2 \le \sqrt{2/\lambda}$. Therefore $\|x^*\|_2 \le \min\left\{\frac{1}{n\lambda}\sum_{i=1}^n\|w_i\|_2, \sqrt{2/\lambda}\right\}$, which finishes the proof. \qed

\bibliographystyle{amsplain}
\bibliography{LF-papers}

\end{document}